\documentstyle[12pt]{article}

\newtheorem{thm}{Theorem}[section]
\newtheorem{defn}[thm]{Definition}
\newtheorem{lemma}[thm]{Lemma}
\newtheorem{prop}[thm]{Proposition}

\newcommand{\beq}[1]{\begin{equation}\label{#1}}
\newcommand{\enq}[0]{\end{equation}}

\newcommand{\qed}[0]{\begin{flushright} \rule{2mm}{3mm} \end{flushright}}
\newcommand{\C}[2]{{{#1}\choose{{#2}}}}
\newcommand{\ga}[0]{\alpha }
\newcommand{\gb}[0]{\beta }
\newcommand{\gc}[0]{\gamma }
\newcommand{\gd}[0]{\delta }
\newcommand{\gD}[0]{\Delta }
\newcommand{\gG}[0]{\Gamma }
\newcommand{\gl}[0]{\lambda }
\newcommand{\gL}[0]{\Lambda}
\newcommand{\go}[0]{\omega}
\newcommand{\gO}[0]{\Omega}

\newcommand{\gs}[0]{\sigma}
\newcommand{\gS}[0]{\Sigma}

\newcommand{\eps}[0]{\varepsilon }

\newcommand{\0}[0]{\emptyset}

\newcommand{\ra}[0]{\rightarrow}
\newcommand{\la}[0]{\leftarrow}
\newcommand{\Ra}[0]{\Rightarrow}

\newcommand{\ZZ}[0]{{\bf Z}}

\newcommand{\II}[0]{{\bf I}}

\newcommand{\Rr}[0]{\mbox{${\bf R}$}}

\newcommand{\Zz}[0]{{\bf Z}}

\newcommand{\g}[0]{{\cal G}}

\newcommand{\I}[0]{{\cal I}}
\newcommand{\J}[0]{{\cal J}}

\renewcommand{\S}[0]{{\cal S}}
\newcommand{\T}[0]{{\cal T}}
\newcommand{\U}[0]{{\cal U}}

\newcommand{\bn}[0]{\bigskip\noindent}
\newcommand{\mn}[0]{\medskip\noindent}
\newcommand{\nin}[0]{\noindent}

\newcommand{\sub}[0]{\subseteq}
\newcommand{\sm}[0]{\setminus}
\renewcommand{\dots}[0]{,\ldots,}

\newcommand{\uzero}[0]{\underline{0}}

\newcommand{\dist}[0]{{\rm dist}}

\newcommand{\gcq}[0]{\gG_Q}
\newcommand{\goj}[0]{G_0^j}
\newcommand{\goo}[0]{G_0'}
\newcommand{\boo}[0]{B_0'}

\newcommand{\sjs}[0]{\gs_j^*}

\newcommand{\indices}[0]{[-d,d]\sm\{0\}}

\newcommand{\gone}[0]{\g^{\star}}
\newcommand{\gtwo}[0]{\g'}
\newcommand{\yx}[0]{2^{\cal E}\times 2^{\cal O}}
\newcommand{\ist}[0]{\I^{\star}}
\newcommand{\sji}[0]{\sigma_j^{-1}}

\newcommand{\sj}[0]{\sigma_j}
\newcommand{\sk}[0]{\sigma_k}

\newcommand{\intb}[0]{\partial^{\star}}

\begin{document}
\renewcommand{\thefootnote}{\fnsymbol{footnote}}

\begin{center}
{\Large\bf
On phase transition in \\
the hard-core model on $\Zz^d$}
\end{center}

\begin{center}
David Galvin\\
Department of Mathematics\\
Rutgers University\\
New Brunswick, NJ 08903\\
and\\
Jeff Kahn\footnotemark\footnotetext{Research
supported in part by
NSF grant DMS-9970433.}\\
Department of Mathematics and RUTCOR\\
Rutgers University\\
New Brunswick, NJ 08903
\end{center}
\date{}

\bigskip\begin{abstract}
It is shown that the hard-core model on ${\bf Z}^d$
exhibits a phase transition at activities above some function
$\lambda(d)$ which tends to zero as $d\rightarrow \infty$; that is:

Consider the usual nearest neighbor graph on ${\bf Z}^d$,
and write ${\cal E}$ and ${\cal O}$ for the sets of even
and odd vertices (defined in the obvious way).
Set
$$\gL_M=\gL_M^d =\{z\in{\bf Z}^d:\|z\|_{\infty}\leq M\},~~~
\intb \gL_M =\{z\in{\bf Z}^d:\|z\|_{\infty}= M\},$$
and write ${\cal I}(\gL_M)$ for the collection of independent
sets (sets of vertices spanning no edges) in $\gL_M$.
For $\lambda>0$ let ${\bf I}$ be chosen from ${\cal I}(\gL_M)$
with $\Pr({\bf I}=I) \propto \lambda^{|I|}$.

\medskip\noindent
{\bf Theorem}
There is a constant $C$ such that if
$\lambda > Cd^{-1/4}\log^{3/4}d$, then
$$\lim_{M\rightarrow\infty}\Pr(\underline{0}\in{\bf I}|{\bf I}\supseteq
\intb \gL_M\cap {\cal E})~>
\lim_{M\rightarrow\infty}\Pr(\underline{0}\in{\bf I}|
{\bf I}\supseteq
\intb \gL_M\cap {\cal O}).$$

\medskip\noindent
Thus, roughly speaking,
the influence of the boundary on behavior at the origin
persists as the boundary recedes.

\end{abstract}

\section{Introduction}

The ``hard-core model" is a simple mathematical model of a gas
with particles of non-negligible size.  The vertices (``sites") of
a graph are regarded as positions, each of which can be occupied
by a particle, subject to the rule that two neighboring sites
cannot both be occupied (particles cannot overlap).

We need a few definitions, but aim to be brief.
For good introductions to the hard-core model see
\cite{BS}, \cite{Haggstrom}.
See also \cite{Georgii} for more general background,
and e.g. \cite{MGT} or \cite{Diestel} for graph theory basics.
A few conventions are mentioned at
the end of this section.

Write $\I(\Sigma)$ for the collection of independent sets (sets of
vertices spanning no edges) of graph $\Sigma$.

For $\gS$ finite and $\gl >0$, the
{\em hard-core measure} with
{\em activity} (or {\em fugacity})
$\gl$ on $\I =\I(\gS)$ (or ``on $\gS$")
is given by
$$
\mbox{$
\mu(I) =
\gl^{|I|}/Z
~~~~$  for  $~ I\in \I$,}
$$
where $Z$ is the appropriate normalizing constant
({\em partition function}),
$Z =\sum\{\gl^{|I'|}:I'\in\I\}$.
(The more usual etiquette here considers probability
measures on $\{0,1\}^{V(\gS)}$ supported on indicators of
independent sets; but the present usage is convenient for
us, and we adhere to it throughout.)

In particular $\gl =1$ gives uniform distribution.
One may also assign different activities
$\gl_v$ to the different vertices $v$ and take $\mu(I)$ proportional to $\prod_{v\in I}\gl_v$,
but we will not do so here; again see \cite{BS},
\cite{Haggstrom}, and also e.g. \cite{K-Kayll}, \cite{lmulti}, \cite{Ded}
for some combinatorial applications.

For infinite $\Sigma$ a measure $\mu$ on $\I(\Sigma)$ is hard-core
with activity $\gl$ if, for $\II$ chosen according to $\mu$ and
for each finite $W\subset V=V(\Sigma)$, the conditional
distribution of $\II\cap W$ given $\II\cap (V\sm W)$ is $\mu$-a.s.
the hard-core measure with activity $\gl$ on the independent sets
of $\{w\in W: w\not\sim \II\cap (V\sm W)\}$ (the vertices that can
still be in $\II$ given $\II\cap (V\sm W)$). General
considerations (see \cite{Georgii}) imply that there is always at
least one such $\mu$; if there is more than one, the model is said
to have a {\em phase transition}.

The canonical (and by far
most studied) case of the hard-core model is
that of (the usual nearest neighbor graph on)
$\ZZ^d$.
Here the seminal result is due to
Dobrushin \cite{Dobrushin}, who proved that there is
a phase transition for
sufficiently large $\gl$, depending on $d$.
(Dobrushin's result was rediscovered by Louth \cite{Louth}
in the context of communications networks.)

The $\gl$ required in \cite{Dobrushin} is
larger than one would expect,\footnote{No explicit bound is
given in \cite{Dobrushin}, but several colleagues report
that Dobrushin's argument works for $\gl>C^d$ for a suitable
constant $C$.} and attempted improvements
have been the subject of considerable effort---if
not publication---in both the
statistical mechanics and discrete mathematics
communities in recent years.

Even the fact that the required $\gl$ increases with $d$
is a little strange, since one expects that as $d$ grows
phase transition should get ``easier," in the sense that
for a given $\gl$, phase transition in dimension $d$
should imply phase transition in all higher dimensions;
but this remains open.

Also open is the existence of a ``critical" activity,
$\gl_c(d)$,
such that one has phase transition
for $\gl >\gl_c(d)$ but not for $\gl < \gl_c(d)$.
While this seems certain to be true for $\Zz^d$,
a cautionary note is sounded in \cite{BHW},
where it is shown that there are graphs (even trees)
for which there is no such critical activity.

As a temporary substitute
we may define
$\gl(d)$ to be the supremum of those $\gl$ for which the hard-core
model with activity $\gl$ on $\ZZ^d$ does not have a phase
transition.

So Dobrushin at least tells us that
$\gl(d) < \infty$, while ``easier as dimension
grows" would imply $\gl(d)<O(1)$.
A particular question that has received much of the attention
devoted to this problem is whether $\gl(d)\leq 1$ for
large $d$.
But in fact it has been generally believed (despite some early
guesses to the contrary) that $\gl(d) $ tends
to zero as $d$ grows; this is what we prove:

\begin{thm}\label{Thm0}
$\gl(d) = O(d^{-1/4}\log^{3/4}d)$.
\end{thm}
The bound here is undoubtedly not best possible;
$O(\log d/d)$ and $O(1/d)$ are
natural guesses at the true value of $\gl(d)$.

\bigskip
We assume henceforth that $d$ is large enough to
support our various assertions.

The problem of showing existence of a phase transition may be
finitized as follows. Let $\gL=\gL_M = \Zz^d\cap [-M,M]^d= {\cal
O}\cup {\cal E}$ with ${\cal O}$ and ${\cal E}$ the sets of odd
and even vertices (defined in the natural way:  $x\in \Zz^d $ is
odd if $\sum x_i$ is odd); let $\mu_M$ be the hard-core measure
with activity $\gl$ on $\gL$ (meaning, of course, on the subgraph
of $\Zz^d$ induced by $\gL$); and (with $\II$ chosen according to
$\mu_M$) let $\mu_M^e$ be $\mu_M$ conditioned on the event
$\{\II\supseteq \intb \gL\cap {\cal E}\}$, where $\intb \gL :=
[-M,M]^d \setminus [-(M-1),M-1]^d$, and define $\mu_M^o$
similarly.

In \cite{BS} it is shown ({\em inter alia}) that the sequences
$\{\mu_M^e\}$ and $\{\mu_M^o\}$ converge to weak limits,
called $\mu^e$ and $\mu^o$, and that there is a phase transition
iff these limits are different.
(This is mainly based on the FKG Inequality, and applies
to general bipartite graphs $\gS$, provided we allow $\{\gL_M\}$ to
be an arbitrary nested sequence with $\cup \gL_M=V(\gS)$.)

Thus it is natural to try to prove phase transition by
exhibiting some
statistic distinguishing $\mu^e$ from $\mu^o$.
We will show
$\mu^e(\uzero\in\II)\neq \mu^o(\uzero\in\II)$, i.e.
\beq{lim}
\lim_{M\ra\infty} \mu_M^e(\uzero\in\II)\neq \lim_{M\ra\infty}
\mu_M^o(\uzero\in\II).
\enq
(Of course we are only using the trivial direction of ``phase
transition iff $\mu^e\neq\mu^o$."
It is not hard to show that (\ref{lim}), too,
is equivalent to phase transition.)

To establish (\ref{lim})
(assuming at least $\gl =\gO(1/d)$,
which is easily seen to be necessary for phase transition)
it is in turn enough to show that for $v_0\in \gL$,
$$
\begin{array}{cccl}
\mu_M^e(v_0) & < & o(1/d) & \mbox{if $v_0$ is odd,} \\
\mu_M^o(v_0) & < & o(1/d) & \mbox{if $v_0$ is even}.
\end{array}
$$

\nin
For then (writing $N$ for neighborhood)

\begin{eqnarray*}
\mu_M^e(\uzero\in \II) &=&
\mu_M^e(N(\uzero)\cap \II=\0)
\mu_M^e(\uzero\in \II|N(\uzero)\cap \II=\0) \\
&=&(1-o(1))\gl/(1+\gl),
\end{eqnarray*}
so that $\mu^e(\uzero\in\II)= (1-o(1))\gl/(1+\gl)$,
whereas $\mu^o(\uzero\in\II)= o(1/d)$.

So in particular the next theorem, whose proof is the main
business of this paper, contains Theorem~\ref{Thm0}.

\begin{thm}\label{Thm}
For
\beq{lambda}
\gl =\go(d^{-1/4}\log^{3/4}d),
\enq
M arbitrary, and $v_0$ an odd vertex of $\gL_M$,
\beq{Pv0}
\mu_M^e(v_0\in\II) < (1+\gl)^{-(2-o(1))d} .
\enq
The same result holds if we reverse the roles of even and odd.
\end{thm}
{\em Remark.}
It is easy to see that

\begin{eqnarray*}
\mu_M^e(v_0\in\II)&=&
\mu_M^e(N(v_0)\cap \II=\0 )\mu_M^e(v_0\in\II|N(v_0)\cap \II=\0)\\
&>&(1+\gl)^{-2d}\frac{\gl}{1+\gl},
\end{eqnarray*}
so that (\ref{Pv0}) actually gives the asymptotics of
$\log \mu_M^e(v_0\in\II)$.

\bigskip

Set
$$
\J =\{I\in\I(\gL):\intb \gL\cap {\cal E}\sub I\}.
$$
The proof of Theorem~\ref{Thm} is a sort of ``Peierls argument"
(see e.g. \cite{Grimmett}):
we try to associate with each $I\in \J$ containing $v_0$
a ``contour"---some kind of membrane separating the outer
even region from an inner odd region containing $v_0$---and
then use this to map $I$ to a large set of
$J$'s, also from $\J$ but not containing $v_0$,
each obtained from $I$ by some modification of the inner region.

This is no surprise:  almost every attempt at settling
this problem that we're aware of
has attacked it more or less
along these lines.  (The one exception is the entropy
approach of \cite{JK}, which for now seems unlikely to
get us to anything like what's proved here.)

The main difficulty in all these attempts has been getting some
kind of control over the set of possible ``contours."
Much of the inspiration for our approach to this problem
was provided by the
beautiful ideas of A. Sapozhenko \cite{Sapozhenko}, which he
used to give, for example, relatively
simple derivations of Korshunov's
\cite{Korshunov} description of the asymptotics for Dedekind's
Problem (in \cite{Sap1}), and, in \cite{Sap2},
of the
asymptotics for the number of independent
sets (``codes of distance 2") in the Hamming cube $\{0,1\}^n$
originally established in \cite{K-S}.

Some of our tools also come from \cite{Sapozhenko}:
Lemma~\ref{App2} is an improved version of one of Sapozhenko's
arguments, and our uses of Lemmas~\ref{Tree}-\ref{Lconn} are
similar to his.

\bigskip
The rest of the paper is devoted to the proof of Theorem~\ref{Thm}.
Unfortunately, saying anything even mildly intelligible about
the argument turns out to be awkward without some preliminaries,
so we will wait:  see the end of Section~\ref{TP} and most of
Section~\ref{Pre}.
(Section~\ref{TP} reformulates slightly and says what we
will actually prove.)

\bn {\em Usage}

We use ``bigraph'' for ``bipartite graph.''

For a graph on vertex set $V$, we use
$\nabla(W)$ for the set of edges having exactly one end in
$W\sub V$ and
$\nabla(U,W)$ for the set of edges having one end in
$U$ and the other in $W$.

The neighborhood of (i.e. set of vertices
adjacent to) $v$ is
$N(v)$; $N(W) =\cup\{N(v):v\in W\}$;
and $\partial W = N(W)\sm W$.
We use $d(\cdot)$ for degree---$d(v)=|N(v)|$ and
$d_W(v)=|N(v)\cap W|$---and $\dist(\cdot,\cdot)$ for distance.

One common abuse:  we often fail to distinguish between a
graph and its set of vertices, so for instance might use
``component" where we should really say
``set of vertices of a component."

When the difference makes no difference, we
pretend that all large numbers are integers.
All constants implied by the notations $O(\cdot)$, $\gO(\cdot)$
are absolute; that is, they do not depend on $d$.

\section{Proof of Theorem~\ref{Thm}}
\label{Proof}

\subsection{Preliminaries}

Here we collect what we will need in the way of known results.

\begin{lemma}\label{Tree}
In any graph with all degrees at most $D$, the number of
connected, induced subgraphs of order $n$ containing a fixed
vertex $x_0$ is at most $(eD)^n$.
\end{lemma}
This follows from the well-known fact (e.g. \cite[p.396,
Ex.11]{Kn}) that the infinite $D$-branching rooted tree contains
precisely $\frac{1}{(D-1)n+1}{ Dn \choose n}$ rooted subtrees of
size $n$.

\bigskip
The next lemma is a special case of a fundamental result due to
Lov\'asz \cite{Lovasz} and Stein \cite{Stein}
(see also \cite{Furedi}).
For a bigraph $\gS$ with bipartition $X\cup Y$, say $Y'\sub Y$
{\em covers} $X$ if each $x\in X$ has a neighbor in $Y'$.

\begin{lemma}\label{Lcor}
If $\gS$ as above satisfies $d(x)\geq a ~\forall x\in X$ and
$d(y)\leq b ~\forall y\in Y$,
then $X$ is covered by some $Y'\sub Y$ of size at most
$(|Y|/a)(1+\ln b)$.
\end{lemma}

Call a set $T$ of vertices of a graph {\em c-clustered} if for any
$x,y\in T$ there are vertices $x=x_0,x_1\dots x_k=y$ with
$\dist(x_{i-1},x_i)\leq c$ for all $i$. The next lemma is from
\cite{Sapozhenko} (see Lemma 2.1); the interested reader should
have no difficulty supplying a proof.

\begin{lemma}\label{Lconn}
If $\gS$ is a graph on V and $S,T\sub V$ satisfy

\mn {\rm (i)}  $S$ is $a$-clustered,

\mn
{\rm (ii)}  $\dist(x,T)\leq b ~\forall x\in S$ and
 $\dist(y,S)\leq b ~\forall y\in T$,

\mn then $T$ is $(a+2b)$-clustered.
\end{lemma}

Finally, we need to know something about isoperimetry in $\ZZ^d$.
Write
$|x|$ for the $\ell_1$-norm of $x$, and set
$B(r)= \{x \in \Zz^d : |x| \leq r\}$,
$S(r)= \{x \in \Zz^d : |x| = r\}$, $b(r)=|B(r)|$ and
$s(r)=|S(r)|$.

\begin{lemma}\label{BL}

Let $C$ be a subset of $\Zz^d$ with

$$|C|=b(r) + \alpha s(r+1),$$
where $0 \leq \alpha < 1$. Then

$$|\partial C| \geq (1-\alpha) s(r+1)+
\alpha s( r+2).$$

\end{lemma}
This is an immediate consequence of
a corresponding
inequality for the torus
$(\Zz / k\Zz)^d $, given by Bollob\'as and Leader
in \cite[Cor. 5]{Bollobas-Leader}.
The case $\ga=0$
was proved by Wang and Wang \cite{WW}.

\subsection{To prove}
\label{TP}

We assume henceforth that $\gl$ satisfies (\ref{lambda}).
We prove only the first part of Theorem~\ref{Thm}
((\ref{Pv0}) for odd $v_0$); switching ``even" and ``odd"
throughout the argument gives the proof of the second part.

It will be convenient to replace the box $\gL_M$ by the discrete
torus $\gG=\gG_M$ obtained from $\gL_M$ by setting
$M=-M$ and identifying
vertices accordingly.
Following our favorite abuse, we regard $\gG$ as either a graph
or a set of vertices as convenient.

We then use $\gD$ for the image of $\intb \gL_M$ under the natural
projection $\gL_M\mapsto \gG$, and continue to write $\uzero$ for
the image of $\uzero$ in $\gG$, and to use ${\cal O}$ and ${\cal
E}$ for the sets of odd and even vertices of $\gG$.

Having done this, we replace $\intb \gL_M$ by $\gD$ in the
definition of $\J$ ($\J=\{I\sub \gG:\mbox{$I$ independent,
$\gD\cap {\cal E}\sub I$}\}$), define $\mu_M^e$, $\mu_M^o$ as
before, and simply regard Theorem~\ref{Thm} as referring to $\gG$,
a change which clearly does not affect its meaning.

\bigskip
We will show a bit more than (\ref{Pv0}): for $I\in \J$, let
$Z=Z(I)$ be the component of $\gG -(I\cap {\cal O})$ containing
$\gD$; then

\beq{pv0}
\mu_M^e(v_0\not\in Z(\II)) <
(1+\gl)^{-(2-o(1))d}
\enq

Let $\J_0 =\{I\in \J: v_0\not\in Z(I)\}$,
and write $w(I)$ for $\gl^{|I|}$.
We prove
(\ref{pv0}) by producing a ``flow"
$\nu:\J_0\times\J\ra [0,1]$ satisfying
\beq{nuout}
\sum_J\nu(I,J)=1  ~~~\forall I\in \J_0
\enq
and
\beq{nuin}
\sum_I\frac{w(I)}{w(J)}\nu(I,J)<
(1+\gl)^{-(2-o(1))d}
~~~\forall J\in \J.
\enq
This gives (\ref{pv0}):
\begin{eqnarray*}
\sum_{I\in \J_0}w(I)
 &=& \sum_{I\in \J_0}w(I)
\sum_{J\in \J} \nu(I,J)\\
&=&
\sum_{J\in \J}w(J) \sum_{I\in \J_0}\frac{w(I)}{w(J)}\nu(I,J)\\
&<&
(1+\gl)^{-(2-o(1))d}
\sum_{J\in \J}w(J).
\end{eqnarray*}

Throughout our discussion we fix $v_0$ and use $I$ for
members of $\J_0$ and $J$ for general members of $\J$.

The definition of $\nu(I,\cdot)$ will depend on a pair
$(G,A)=(G(I),A(I))\in 2^{\cal E}\times 2^{\cal O}$ associated with
$I$. The construction and salient properties of the pair are given
in Sections~\ref{Con} and \ref{SV}, but it will not be until
Section~\ref{Flow} that we are able to specify $\nu$. First steps
toward this specification are taken in Section~\ref{Shifts}, which
finally puts us in a position---in Section~\ref{Pre}---to give
some clue as to how the main part of the argument will proceed.

\subsection{``Contours"}
\label{Con}

For a set $P$ of vertices (in any graph)
we use $\intb P$ for the
{\it internal boundary} of $P$:
$$\intb P = \{v \in P | N(v) \not \sub P \}.$$

The following observation is used several times, so
we record it as a lemma; its easy proof
is left to the reader.

\begin{lemma}\label{bdL}
Let $\gS$ be a graph, $S\sub V(\gS)$, and
$T$
(the vertex set of)
some component of $\gS-(S\sm \intb S)$.
Then $\intb T \sub \intb S$.
\end{lemma}

Let $I\in\J_0$, $Z=Z(I)$ be as in Section~\ref{TP}, and set $Z_0 =
\intb Z$. By the definition of $Z$, it is clear that $Z_0\subset
{\cal E}$ and $Z_0\cap I=\0$. Let $W'$ be the component of $v_0$
in the graph $\Gamma - (Z \setminus Z_0)$. By Lemma~\ref{bdL},
$\intb W' \sub W'\cap Z_0\sub {\cal E}$.

Let $W'' = W' \cup \{x \in {\cal O} | N(x) \sub W' \}$. This is
clearly connected, with $ \intb W''\sub \intb W'$.

Now consider $\Gamma - (W'' \setminus \intb W'')$. This breaks
into a number of components, one of which, $C$ say, contains
$\gD$. Again using Lemma~\ref{bdL}, we have $\intb C \sub C\cap
\intb W''. $ Finally, set $W = \gG\sm (C\setminus \intb C)$, $G=W
\cap {\cal E}$, $A=W \cap {\cal O}$, and $G_0 = \intb W$.

The next proposition collects relevant properties of these
objects.
Once we have these properties, we will not
be concerned
with how $G,A$ etc. were derived from $I$.

\begin{prop}\label{GAproperties}

\beq{GA0}
v_0\in A; ~~~ W\cap \gD=\0;
\enq

\beq{GA1}
\mbox{both $C$ and $W$ are connected;}
\enq

\beq{GA2}
G_0 = \intb C;
\enq

\beq{GA3} \mbox{$G = N(A)~~$ and $~~ A = \{x \in {\cal O} | N(x)
\sub G\}$;} \enq

\beq{GA5}
G_0 \cap I = \emptyset;
\enq

\beq{GA6}
N(G_0) \cap I \subset A;
\enq

\beq{GA7}
G_0\sub N(A\cap I).
\enq

\end{prop}
{\em Proof.}
Both (\ref{GA0}) and the connectivity of $C$ are immediate.
To see that $W$ is
connected, notice that each component of $\gG-(W''\sm \intb W'')$
must meet $\intb W''$ (or it would be a component of the
connected graph $\gG$).
Thus $W$ is the union of the connected set $W''$ and a number
of other connected sets each of which meets $W''$, so is itself
connected.  So we have (\ref{GA1}).

For (\ref{GA2}): $\intb C \sub W\cap {\cal E}$ and the
connectivity of $C$ give
$$x \in \intb C\Ra
\0\neq N(x)\cap C\sub C\cap {\cal O}\sub C\sm W \Ra x\in \intb
W,$$ so $\intb C\sub\intb W$; and Lemma \ref{bdL} and the
connectivity of $W$ give the reverse containment.

Connectivity of $W$ and the fact that $G_0\sub {\cal E}$ give
$G=N(A)$. That $ A \sub \{x \in {\cal O} | N(x) \sub G\} $ follows
from $G=N(A)$ (or just $\intb W\sub {\cal E}$). For the reverse
containment, notice that $x\not\in W\Ra N(x)\cap W\sub G_0\sub
W'$, whereas $N(x)\sub W'$ would imply $x\in W''\sub W$; so
$x\not\in W\Ra N(x)\not\sub W$.

For (\ref{GA5}) recall that
$G_0=\intb C\sub \intb W''\sub \intb W'\sub Z_0$ and $Z_0\cap I=\0$.

That $N(G_0) \cap I \sub A$ follows from $G_0\sub \intb W'$,
since $N(\intb W')\cap  I $ is clearly contained in $A$.

Finally, $v \in G_0 \Rightarrow v \in Z_0 \Rightarrow v \sim I$,
so (\ref{GA7}) follows from
(\ref{GA6}).\qed

\subsection{Topology}
\label{SV}

The purpose of this section is to prove, for any $I\in\J_0$ and
$W$, $G$ etc. produced from $I$ as in Section~\ref{Con},

\beq{G0} \mbox{$G_0$ is 2-clustered} \enq Our proof of this, which
is considerably longer than we would wish and unrelated to the
methods in the rest of the paper, might profitably be skipped on a
first reading.

Though (\ref{G0}) turns out to follow from the
connectivity of $W$ and $C$ (see (\ref{GA1})), we could not see a simple
combinatorial proof of the implication, and our argument
requires a little topological detour, based on

\begin{lemma}\label{AlgTop}
If $U,V$ are connected subsets of $X=\Rr^n$ or $S^n$, $n>1$,
with $U\cup V=X$, $U$ closed and $V$ compact, then
$U\cap V$ is connected.
\end{lemma}
(As usual, $S^n$ is the unit sphere $\{x\in\Rr^{n+1}:\sum x_i^2=1\}$.
We also write $B^{n+1}$ for the corresponding unit ball.)

The (presumably well-known) proof of Lemma~\ref{AlgTop}
is given at the end of this section.

It will  be convenient here to write
$\Omega$ for the nearest neighbor graph on $\Zz^d$.
As usual, $\Omega[S]$ is the subgraph induced by $S$.
We will prove (\ref{G0}) in the following more general form.

\begin{prop}\label{connbound}

Let $R\cup B$ be a decomposition of $V(\Omega)$ ($=\Zz^d$), with
both $\Omega[R]$ and $\Omega[B]$ connected and $R$ finite. Suppose
$G:=R\cap B$ is contained in ${\cal E}$ and is the internal
boundary of each of $R,B$.  Then $G$ is 2-clustered.

\end{prop}
{\em Remark.} We will actually show that $G$ is $2$-clustered in
each of $R$ and $B$.

\mn {\em Proof} With $\Omega$ embedded in $\Rr^d$ in the natural
way, we extend $R$ and $B$ to closed connected subsets $R^*$ and
$B^*$ of $\Rr^d$ so that $R^* \cup B^* = \Rr^d$ and $G^* := R^*
\cap B^*$ is path-connected. We then derive the $2$-clusteredness
of $G$ from the path-connectedness of $G^*$.

We view $\Rr^d$ as the union of
$\Zz^d$-translates of $[0,1]^d$ (the {\em cells} of $\Rr^d$),
and define $R^*$ and $B^*$ cell by cell.
Within a cell we proceed by dimension, first defining the
extensions for
$0$-dimensional faces (the
vertices of $\Omega$), $1$-dimensional faces (the edges of
$\Omega$), and $2$-dimensional faces, and then
continuing inductively.
(As usual a {\em face} of a cell is the intersection of the
cell with some supporting hyperplane.
Henceforth we use ``$k$-face" for ``$k$-dimensional face.")
For the inductive
step, we need a topological lemma (Lemma~\ref{extending}),
for the statement of which it's convenient to
introduce two local definitions.
Let us say that a subset of a
topological space is {\em civilized} if it is
closed, has only finitely many components, and each of its
components is path-connected.

\begin{defn}\label{nice}

A decomposition $X = R \cup B$ of a topological space $X$, with
$R\cap B=G$, is {\em nice} if it satisfies:

\mn
{\rm (i)} $G=\partial R = \partial B$;

\mn
{\rm (ii)}  each of $R$, $B$, $G$ is civilized; and

\mn
{\rm (iii)}  each of $R$, $B$---and so each component of
$R$ and $B$---is the closure of the
union of finitely many open, path-connected sets.

\mn
If $X=R\cup B$ is a nice decomposition, and $R'$, $B'$
are obtained from $R$, $B$ by adding finitely many points,
then we also call the decomposition $X=R'\cup B'$ nice.

\end{defn}
(Of course there is some redundancy in conditions (i)-(iii).)

We say that two nice decompositions $X_1 = R_1 \cup B_1$ and $X_2 =
R_2 \cup B_2$ are {\em compatible} if $R_1\cap X_1 \cap X_2 =
R_2\cap X_1 \cap X_2$ and $B_1\cap X_1 \cap X_2 = B_2\cap X_1 \cap X_2$.
It's straightforward
to check that nice decompositions of different spaces
can be combined if they are compatible:

\begin{lemma}\label{combining}

Suppose $X=X_1\cup\cdots\cup X_m$ with each $X_i$ closed.
If $X_i=R_i\cup B_i$ are pairwise compatible, nice decompositions,
then $(\cup R_i)\cup (\cup B_i)$ is a nice decomposition of $X$.
\end{lemma}

We now state the topological lemma alluded to above,
deferring its proof until after the derivation of
Proposition~\ref{connbound}.
(Recall
$B^{n+1}$ and
$S^n$ are the unit ball and sphere in $\Rr^{n+1}$.)

\begin{lemma}\label{extending}

Assume $n>1$.
If $R \cup B$ is a nice decomposition of $S^n$, then there is
a nice decomposition $ R^* \cup B^*$ of  $B^{n+1}$,
with $R^* \cap S^n = R$, $B^* \cap S^n = B$,
and such that if $C$ is any
component of $R^*$ (resp. $B^*$, $G^*$), then $C \cap S^n$ is a
component of $R$ (resp. $B$, $G$).

\end{lemma}
(This is easily seen to fail for $n=1$.
It may be worth pointing out that for $R$ and $B$, condition
(iii) of Definition~\ref{nice} refers to sets that are open
{\em in} $S^n$; similarly $\partial R$ and $\partial B$ are
boundaries relative to $S^n$, while $\partial R^*$ and $\partial B^*$
are boundaries relative to $B^{n+1}$.)

Of course Lemma~\ref{extending} still applies if we replace
the $B^{n+1}$ by any of its homeomorphic images
(and $S^n$ by the corresponding homeomorphic copy);
in our case the relevant image will be $[0,1]^d$.

\medskip
We now fix a cell, and begin defining our extensions.
For vertices and edges we do the natural things:
$R^*\cap V(\Omega) =R$, $B^*\cap V(\Omega) =B$;
and we put (the interior of) an edge in $R^*$ (resp. $B^*$)
iff both its ends are in $R^*$ (resp. $B^*$),
noting that exactly one of these possibilities
occurs, since $\nabla(G,G)=\emptyset$.

Next, we deal with $2$-dimensional faces. If the vertices of such a
face are all in $R$ (resp. $B$), then put the interior of the
face in $R^*$ (resp. $B^*$). Otherwise, the face has two
opposite corner vertices ($v_1, v_3$, say) in $G$,
with one of its remaining
two vertices ($v_2$) in $R \setminus B$ and the other
($v_4$) in
$B \setminus R$. Put the interior of the convex hull of $v_1, v_2,
v_3$ in $R^*$, the interior of the convex hull of $v_1, v_3, v_4$ in
$B^*$, and the interior of the diagonal joining $v_1$ and $v_3$ in
$R^* \cap B^*$. It is easy to check that these
$(R^*, B^*)$-decompositions of the
$2$-dimensional faces are nice.
(It may be worth observing that a 2-dimensional face contained
in $R^*$ may still have one or two of its vertices in $B^*$,
and vice versa.)

We now proceed by induction, assuming the decomposition
has been defined on faces of dimension less than
$k\in\{3\dots d\}$.
Each $k$-face $F$ is homeomorphic to $B^k$, and is bounded by
the union of finitely many $(k-1)$-dimensional faces. The
decomposition of each of
these bounding faces is nice, and the decompositions on any two faces
are compatible (since we are defining the decomposition from lower
dimensions up). So, by Lemma~\ref{combining}, we have a nice
decomposition of the boundary of $F$.  We now
apply Lemma \ref{extending} to extend to a nice decomposition of the
entire face. Once we have a nice decomposition of each cell, we
get the full decomposition $\Rr^d=R^*\cup B^*$ by
combining the
decompositions of the cells, again appealing to
Lemma~\ref{combining} for ``nice."
(For formal applicability of the lemma, we can use a single $X_i=B_i$
for the union of all cells not meeting $R$.)

It is clear from the construction that $R^*$ and $B^*$ are closed,
$R^*$ is bounded, and $R^* \cup B^* = \Rr^d$. To see that $R^*$ is
connected, notice that by
construction, any component of $R^*$ contains an edge of $\Omega[R]$,
and that every edge of $\Omega[R]$ is contained in a component of
$R^*$; connectivity of $R^*$ then follows from connectivity of
$\Omega[R]$. The same argument shows that $B^*$ is
connected.

Lemma~\ref{AlgTop} now shows that $G^*$ is connected,
which, since $G^*$ is also civilized (since $R^*\cup B^*$ is nice),
implies that it is actually path-connected.

It remains to show that path-connectedness of $G^*$ implies
2-clusteredness of $G$. It is enough to show that for each pair of
vertices $u, v \in G$, there is a path connecting them in $G^*$
which is supported entirely on the $2$-dimensional faces of
$\Rr^d$; for, by the construction of $R^*$ and $B^*$, such a path
is supported on diagonals (of $2$-dimensional faces) connecting
pairs of vertices from $G$, and such diagonals correspond to steps
of length $2$ in $\Omega$. (This also justifies the remark
following Proposition \ref{connbound}.)

So, consider a $(u,v)$-path $P$ in $G^*$ given by the
continuous function $f:[0,1] \rightarrow \Rr^d$.
If $P$ is supported
on $2$-dimensional faces of $\Rr^d$, then we are done.
Otherwise, let $k>2$ be the maximum dimension of a face
whose interior meets $P$.
It's enough to show that we can replace $P$ by a path
meeting
the interiors of fewer $k$-faces than $P$
and no faces of dimension more than $k$.

To do this,
choose a $k$-face $F$ and component $C$ of $G^*\cap F$ with
$C\cap F^0\cap P\neq \0$
(where $F^0$ is the interior of $F$).
Let $p=\inf \{x \in
[0,1]:f(x) \in C\cap F^0 \}$ and
$q = \sup \{x \in [0,1]:f(x) \in C\cap F^0 \}$.
Then $f(p), f(q) \in C \cap \partial F$, which, by
construction, is path-connected. So we may replace $f([p,q])$
in $P$ by a path contained in $\partial F$.\qed
{\em Proof of Lemma \ref{extending}}

To avoid confusion, we now write $\partial X$,
$\partial\hspace{.01in}' X$ and
$\partial\hspace{.01in}'' X$ for the boundaries of $X$ relative to,
respectively, $\Rr^{n+1}$, $B^{n+1}$ and $S^n$.

We may assume neither $R$ nor $B$ contains isolated
points:
otherwise we can simply delete such points, produce $R^*$ and $B^*$
for the resulting ``reduced" $R$ and $B$, and then add the
deleted points of $R$ ($B$) to $R^*$ ($B^*$).

We use $(R,B)$-{\em component} to mean a component of either $R$
or $B$,
and proceed by induction on the number of ($R$,$B$)-components
in the decomposition of $S^n$.

If there is exactly one such component (a
component of $R$, say), then $R=S^n$, and $B=\0$.
Setting $R^*=B^{n+1}$ and $B^*=\emptyset$, we
get a nice decomposition of $B^{n+1}$ which satisfies the conditions
of the lemma.

Otherwise, there must be at
least one ($R$,$B$)-component $T$ for which $S^n \setminus T^0$
is connected. For suppose $S^n \setminus T^0$ is disconnected for
every ($R$,$B$)-component $T$. Choose an ($R$,$B$)-component
$T_0$ ($\sub  R$, say)
such that one of the components of $S^n \setminus
T_0^0$, $C$ say, contains as few $(R,B)$-components as possible,
and let $T_1$ be an $(R,B)$-component of $C$ (i.e. contained
in $C$, noting that each $(R,B)$-component other than $T_0$
is either contained in or disjoint from $C$).
Now $S^n \setminus C^0$ is connected in
$S^n \setminus T_1^0$,
so $S^n \setminus T_1^0$ (which by assumption is not connected)
contains a component whose
($R$,$B$)-components form a proper subset of the
($R$,$B$)-components of $C$, contradicting the choice of $T_0$.

Let $T$, then, be an ($R$,$B$)-component with $S^n \setminus T^0$
connected. We may assume that $T$ is a component of $R$.
Applying Lemma~\ref{AlgTop} with $X=S^n$, $U=T$ and $V=S^n\sm T^0$,
we find that $\partial\hspace{.01in}''T$
is connected, so that $T$ meets
exactly one component, say $C$, of $B$
(and $C\supseteq \partial\hspace{.01in}''T$).

Set $T^*=\{\gl x:x\in T, \gl \in [1/2,1]\}$.
This will be one component of $R^*$.
It is easy to see that $T^*$ is closed and path-connected
(so civilized),
as is $\partial'T^*$,
and that $T^*\cap S^n=T$, a component of $R$.

Now let $(T^*)^0$ be the relative interior of $T^*$ with respect
to $B^{n+1}$ (namely,
$(T^*)^0 = \{ \gl x : x \in T^0,  \gl \in (1/2,1] \}$),
$P=\partial(B^{n+1}\sm (T^*)^0) $
($= (S^n\sm T^0)\cup \partial\hspace{.01in}'T^*$), and
$Q=B^{n+1}\sm (T^*)^0$.
Then $(Q,P)$ is (easily seen to be)
homeomorphic to $(B^{n+1}, S^n)$.

Let, further,
$R_1=R\sm T$, $B_1=B\cup \partial\hspace{.01in}'T^*$,
and $C_1=C\cup \partial\hspace{.01in}'T^*$.
Then

\mn
(i)  the components of $R_1$ are precisely the components of $R$
other than $T$,

\mn
(ii)  the components of $B_1$ are $C_1$ and the components of $B$ other
than $C$,

\mn
and it is easy (if tedious) to deduce that
$R_1\cup B_1$ is a nice decomposition of $P$.

Our inductive hypothesis thus gives a nice decomposition
$R_1^*\cup B_1^*$ of $Q$, and we obtain the desired decomposition,
$R^*\cup B^*$,
of $B^{n+1}$ by setting $B^*=B_1$ and $R^* =R_1\cup T^*$
(again an easy verification using (i) and (ii)).\qed

\bn
{\em Proof of Lemma~\ref{AlgTop}}

\medskip
We first establish a corresponding statement for open sets:
{\em if $U,V$ are connected, open subsets of $X=\Rr^n$ or $S^n$, $n>1$,
with $U\cup V=X$,
then $U\cap V$ is connected.}

\mn
{\em Proof}.
We use the Mayer-Vietoris sequence. If $X$ is a
topological space, and $U$ and $V$ are open subsets of $X$ whose union
is $X$, then this is a long exact sequence of group homomorphisms
ending with

$$\cdots \rightarrow H_1(X) \rightarrow H_0(U \cap V) \rightarrow
H_0(U) \oplus H_0(V)
\rightarrow H_0(X) \rightarrow 0,$$

\noindent where $H_m$ is the $m^{th}$ homology group. We apply
this with $X=\Rr^n$ or $S^n$. Using the facts that $H_m(\Rr^n) =
0$ whenever $m \geq 1$ and that if $O$ is an open subset of
$\Rr^n$ or $S^n$, then $H_0(O) \cong \Zz$ iff $O$ is connected,
this long exact sequence becomes

$$0 \rightarrow H_0(U \cap V) \rightarrow \Zz \oplus \Zz \rightarrow
\Zz \rightarrow 0.$$

From the exactness of this sequence, it follows that $H_0(U \cap V)
\cong \Zz$, so that $U \cap V$ is connected.\qed

Now let $U,V$ be as in the lemma, and for each $\eps >0$, set
$U_\eps = \{x \in X : d(x,U)<\eps \}$ and $V_\eps = \{x \in X :
d(x,V)<\eps \}$. These are open, connected sets whose union is
$X$, so by the preceding result, $U_\eps \cap V_\eps$ is
connected. Thus $\overline{U_\eps \cap V_\eps}$ is connected; it
is also closed and bounded, so compact. So $U \cap V = \cap_{\eps
> 0} \overline{U_\eps \cap V_\eps}$ is the intersection of a
nested sequence of compact, connected sets, so is itself
connected.\qed

\subsection{Shifts and $\varphi_j$}
\label{Shifts}

We again fix $I\in\J_0$ and take $W,G,A$ etc. to be
as in Section~\ref{Con}.

For $j\in
\{\pm 1\dots \pm d\}$, define $\sj$, the
{\it shift in direction $j$},
by

$$\sj (v) = v +  e_j ,$$
where $e_j$ is the $j^{th}$ standard basis vector
if $j>0$ and $e_j=-e_{-j}$ if $j<0$, and
set
$$\goj = \{v \in G_0 : \sj^{-1}(v) \not \in A \}
=G_0\cap \sj ({\cal O}\sm A).$$

\begin{prop}\label{shift}

For each j, the sets
$ I \setminus W$,  $\sj(I\cap W) $ and $\goj $
are pairwise disjoint, and their union is an
independent set.
\end{prop}
{\em Proof.}
Trivially, $\sj(I)\cap I =\0$, so in particular
$(I\sm W)\cap \sj(I\cap W)=\0$;
$(I\sm W)\cap \goj =\0$ is trivial (because $\goj\sub W$);
and $\sj(I\cap W)\cap \goj=\0$
follows from the definiton of $\goj$.
So the union is disjoint.

Clearly $(I \setminus W)$, $\sj(I \cap W)$ and $\goj$ are
all independent sets. To show independence of the union,
we must show that there are no edges between any two of them.
Since $\nabla(I\sm W,W)=\0$ (by (\ref{GA6}))
and $\sj(I\cap W)\sub W$ (by (\ref{GA5})), we have
$\nabla((I\sm W), (\sj(I\cap W)\cup \goj))=\0$.

This leaves $\nabla ( \sj(I \cap W),\goj)$. Suppose, for a
contradiction, that $y\in \goj$ and $\sk(y) \in \sj(I \cap W)$ for
some $k$. Then $z:=\sj^{-1}(\sk(y)) \in I \cap W \cap {\cal
E}\subset G \setminus G_0$ (by (\ref{GA5})), implying
$\sj^{-1}(y)=\gs_k^{-1}(z)\in A$, contrary to the assumption $y\in
\goj$.  So $\nabla (\sj(I \cap W),\goj) = \emptyset$. \qed

Define
$
\sjs(I)= (I\sm W)  \cup \sj(I\cap W)
$
and
$$
\varphi_j(I) = \{J: \sjs(I)\sub J\sub \sjs(I)\cup \goj\}.
$$
Then Proposition~\ref{shift} implies
$$
\varphi_j(I)\sub \J.
$$
Notice also that
we recover $I$ from $j$, $J$ ($\in \varphi_j(I)$) and $(G,A)$;
namely, if we are given $(G,A)$, $j$, and
$J\in\varphi_j(I)$, then
\beq{recover1}
I=(J\sm W)\cup \gs_j^{-1}(J\cap (W\sm G_0^j)).
\enq

\subsection{Conventions and preview}\label{Pre}

{\em Conventions}

In much of what remains we can ignore $I$ and
concentrate on pairs from
$$\g:=\{(G,A)\in\yx:\mbox{$(G,A)$ satisfies (\ref{GA3}) }\}.$$
Notice that under (\ref{GA3}) each of $G$, $A$ determines
the other.

If $(G,A)$ is produced from $I$ as in Section~\ref{Con}
then we write $(G(I),A(I))$, noting that a given
$(G,A)$ may correspond to more than one $I$.

We will always take $W=G\cup A$ and $G_0=\intb W$ (a subset of
${\cal E}$ because of (\ref{GA3})).

Set $\ell = 2d$; so $\gG$ is an $\ell$-regular bigraph.
(We tend to think in terms of $d$ and
use $\ell $ sparingly, for instance usually
preferring $O(d)$ to the equivalent $O(\ell)$.)
Though we usually work in $\gG$, we sometimes---especially
in Section~\ref{SA}---consider more general graphs $\gS$,
always assumed to satisfy

\beq{graph} \mbox{$\gS$ is an $\ell$-regular bigraph with
bipartition $V={\cal O}\cup {\cal E}$.} \enq

We always take $|G|=g$ and $|A|=a=(1-\gd)g$,
and for given $g,\gd$ set
$$\g(g,\gd) =\{(G,A)\in \g:|G|=g,|A|=(1-\gd)g\},$$
$$\J(g,\gd)=\{I\in \J_0:(G(I),A(I))\in\g(g,\gd)\}.$$
(It's generally best to think of $\gd$ as small,
though it will not always be so.)

As will appear, the quantity that really matters is
almost always $\gd g$ ($=|G|-|A|$),
and it will be convenient to take, for any $t$,
$$
\g(t)= \{(G,A)\in\g: |G|-|A|=t\}.
$$
Notice that for $(G,A)\in \g(t)$, \beq{nabla} |\nabla(W,V\sm W)|
~~( = |\nabla(G_0,{\cal O}\sm A)|) ~~ =t\ell. \enq

Though we don't really need $t$, we use it to emphasize a certain
{\em duality}: if $(G,A)\in\g(t)$ in some graph $\gS$ satisfying
(\ref{graph}), then $({\cal O}\sm A,{\cal E}\sm G)$ belongs to the
analogue of $\g(t)$ obtained by reversing the roles of ${\cal O}$
and ${\cal E}$ in $\gS$---but of course $g$ and $\gd$, unlike $t$,
are not usually preserved by this switch.

\bn
{\em Preview}

Our tasks are to define $\nu$,
for which (\ref{nuout}) will turn out to be obvious,
and establish (\ref{nuin}).

We will eventually associate with each $(G,A)$ a particular
index $j=j(G,A)$, and set $j(I)=j(G(I),A(I))$.
(This is basically a $j$ for which $|\goj|=\log_2|\varphi_j(I)|$
is large, though there are some
additional considerations.)
We then define $\varphi(I)=\varphi_{_{j(I)}}(I)$ and require
\beq{phi1}
J\not\in\varphi(I)\Rightarrow \nu(I,J)=0.
\enq

Let us call $I$ {\em small} if $|G(I)|\leq d^3$
(we could get by with $d^{9/4}$;
see (\ref{note3})), and {\em large} otherwise.

For small $I$---an easy case, as we will see in
Section~\ref{Fin}---we simply choose $j=j(I)$ to maximize $|\goj|$
(where $G=G(I)$), so that, since \beq{Gsum} \sum_j |\goj|
=|\nabla(G,{\cal O}\sm A)| =\gd g \ell, \enq we have \beq{gojbig}
|\goj|\geq \gd g. \enq We then set \beq{easynu} \nu(I,J) =
\gl^{|J|-|I|}(1+\gl)^{-|\goj|} ~~~\forall J\in\varphi(I). \enq
(Note this satisfies (\ref{nuout}). The separate treatment of
small $I$ is unnecessary if we only want the phase transition, but
is needed for the ``correct" bound in (\ref{Pv0}).)

Most of our work (including
everything in Sections~\ref{SV} and \ref{CB}-\ref{POM})
is geared to
large $I$ (though often valid in general).
For most of our discussion we fix
$(g,\gd)$, and aim to bound the contribution
of $\J(g,\gd)$ to (\ref{nuin}).
Of course these contributions must eventually be summed,
but this turns out not to add anything significant.

Before beginning in earnest, we pause in
Section~\ref{Iso} to adapt the isoperimetric Lemma~\ref{BL}
to our situation
(Lemma~\ref{deltag}).  This is needed especially in Section~\ref{Fin},
but will also make an appearance in Section~\ref{CB}.

In Sections~\ref{CB}-\ref{Status}
we associate with
each relevant $(G,A)$
some $(F,S)\in \yx$ which ``approximates" $(G,A)$
in an appropriate sense.
The definitions of
$j(I)$ and $\nu(I,\cdot)$
(in Section~\ref{Flow}) are then based on our approximation
to $(G(I),A(I))$.
The main points are:  (i)
the set
of possible approximations is small
(Lemma~\ref{LApp3}); and (ii)
for a given $J$, $I$'s for which $(G(I), A(I))$ is approximated
by a particular
$(F,S)$
don't contribute too much in (\ref{nuin})
(see (\ref{MAIN})), construction of a $\nu$ achieving
this being made
possible by the accuracy of our approximations.

The proof that $\nu$ behaves as desired (that is, of
(\ref{MAIN})) is given in
Section~\ref{POM}, and Section~\ref{Fin} is a mopping up
operation, combining what we already know for large $I$'s with the
easy analysis for small $I$'s and the isoperimetric information
from Lemma~\ref{deltag}, to finally establish (\ref{nuin}).

\bn
{\em More conventions}

For whatever $G,A,F,S$ we have under discussion, we set $H={\cal
E}\sm G$, $B={\cal O}\sm A$, $E={\cal E}\sm F$, $T={\cal O}\sm S$,
$B_0= B\cap N(G)$, $S_0=S\cap N(E)$, and $E_0 = E\cap N(S)$.

\medskip
{\em From now until Section~\ref{Fin}
we fix $g,\gd$ and always take $I \in\J(g,\gd)$ and
$(G,A)\in \g(g,\gd)$.}
(We will not see $I$ again until Section~\ref{Flow}.)

\subsection{Isoperimetry}\label{Iso}

Before continuing,
we need to work out what Lemma~\ref{BL} implies in the
way of a lower bound on $\gd$ for given $g$.

\begin{lemma}\label{deltag}

Suppose
$(G,A)\in\g(g,\gd)$ satisfies
\beq{inside}
(G\cup A)\cap \gD=\0.
\enq
Then

$$
\delta = \left\{ \begin{array}{ll}
                   \gO(g^{-1/d}/d) & \mbox{ for all $g$}\\
                   1-O(1/d) & \mbox{ if $g < d^{O(1)}$.}
                 \end{array}
         \right.
$$

\end{lemma}
(For the $(G,A)$'s of interest to us, (\ref{inside}) is given by
(\ref{GA0}).)

\mn
{\em Proof.}  In view of (\ref{inside}), the lemma does not
change if we replace the torus $\gG$ by the box $\gL$.

For the first part of the
lemma, the main thing we have to show is
\begin{prop}\label{goal1}
$~~~~s(r) =\gO (b(r)^{1-1/d})$
\end{prop}
(where $B(r)$, $S(r)$, $b(r)$, $s(r)$ are as
defined before Lemma~\ref{BL}).
Notice that this, combined with
Lemma \ref{BL}, implies
that for any $C \subset \Zz^d$,

\beq{BL3}
|\partial C| =\gO(|C|^{(d-1)/d}).
\enq

Proposition~\ref{goal1}
is again something for which one would hope
to just give a reference; but we could not find one,
or even give the short proof that seems called for.

For the proof,
we'll be interested in the average number of nonzero
entries in an element of $S(q)$,
$$t(q) := s(q)^{-1}\sum_{x\in S(q)}|supp (x)|.$$
This is useful because, setting
$$N(q)=|\{(x,y)\in S(q)\times S(q+1):x\sim y\}|,$$
we have
$$s(q)(2d-t(q))=N(q)\leq s(q+1)\min\{q+1,d\},$$
implying
\beq{qratio}
\frac{s(q)}{s(q+1)}\leq\frac{\min\{q+1,d\}}{2d-t(q)}.
\enq

This already implies Proposition~\ref{goal1} for, say,
$r\leq .9d$,
since in this case we have
$$
b(r)\leq s(r)\sum_{i=0}^r\frac{(r)_i}{(2d-r+i)_i}
\leq s(r) \sum_{i\geq 0}\left(\frac{r}{2d-r}\right)^i =O(s(r)).
$$

For larger $r$ we will have to work harder.
Here we first show, for $q=\gb d $ with $\gb >.9$,
\beq{t(q)}
t(q) < (1-1/(20\gb))d.
\enq

Let
$$S(q,t)=\{x\in S(q):|supp(x)|=t\},$$
$s(q,t)=|S(q,t)|$,
and define $B(q,t)$ and
$b(q,t) $
similarly.  Then
$$
f(q, t):=\frac{s(q,t+1)}{s(q,t)}= 2\frac{(d-t)(q-t)}{(t+1)t}.
$$
Set
$t_0=t_0(q)=\lceil (1-1/(4\gb))d\rceil$.
Then $t\geq t_0$ implies
\begin{eqnarray*}
f(q,t)&\leq &
2\frac{(1/(4\gb))(\gb -1 +1/(4\gb))}{(1-1/(4\gb))^2}\\
&=& 2\left(\frac{2\gb-1}{4\gb-1}\right)^2 <\frac{1}{2}.
\end{eqnarray*}
Thus
\begin{eqnarray*}
t(q)&=&s(q)^{-1}\sum_{t\leq q}ts(q,t)\\
&<& t_0 +\sum_{i\geq 1}i2^{-i} = t_0+2.
\end{eqnarray*}
This gives (\ref{t(q)}) provided $\gb\leq d/15$.
For larger $\gb $ we just use
$$\frac{s(q,d-1)}{s(q,d)} =\frac{d(d-1)}{2(\gb d-d+1)}>
\frac{d-1}{2\gb},$$
whence
\begin{eqnarray*}
d-t(q) &=& s(q)^{-1}\sum (d-i)s(q,i)
\geq
\sum_{i<d}s(q,i)/(\sum_{i\leq d}s(q,i))\\
&\geq & s(q,d-1)/(s(q,d-1)+s(q,d))
\geq (d-1)/(2\gb +d-1),
\end{eqnarray*}
which again gives (\ref{t(q)}).

Now let $r=\gc d \geq .9d$.  By
(\ref{t(q)}) and (\ref{qratio}) we have, for
$r-i \geq .9d$,
$$
s(r-i)\leq s(r)\prod_{j=1}^i\frac{d}{d+d^2/(20(r-j))}
< s(r)(1-\gO(1/\gc))^i,
$$
so
\beq{b(r)}
b(r)\leq  s(r)\sum_{i=0}^{r-.9d} (1-\gO(1/\gc))^i
 + b(.9d)
= O(\gc s(r))
\enq
(since we know $b(.9d) =O(s(.9d))=O(s(r))$).

On the other hand, with $t_0=t_0(r)$, we have
$$
b(r) > b(r,t_0) =2^{t_0}\C{d}{t_0}\C{r}{t_0}
> \exp [t_0\log(r/t_0)],
$$
and
$b(r)^{1/d}>\exp[(1-1/(4\gc))\log (r/t_0)]=\gO(\gc)$;
and this with (\ref{b(r)}) gives Proposition~\ref{goal1}.\qed

Now for the first part of Lemma~\ref{deltag}, we consider the
possibilities $|G_0| > |A|$ and $|G_0| \leq |A|$ separately, in
both cases using the fact that $|G_0| \leq \delta gd$ (since
$|G_0|\leq |\nabla(G,{\cal O}\sm A)|=\gd g d$).

If $|G_0| > |A|$, then $\delta > 1/(d+1)$, so certainly $\delta =
\gO(g^{-1/d}/d)$. If, on the other hand, $|G_0| \leq |A|$, then we
have (using (\ref{BL3}) and the fact that $\partial ((G \setminus
G_0) \cup A) = G_0$)

\begin{eqnarray*}
\delta & \geq & |G_0|/(dg) \\
 & = & \gO(|(G \setminus G_0) \cup A|^{(d-1)/d}/(dg)) \\
 & = & \gO(|G|^{(d-1)/d}/(dg)) \\
 & = & \gO(g^{-1/d}/d).
\end{eqnarray*}

For small $g$ notice that for $r < O(1)$,
$$s(r) = 2^rd^r/r! + O(d^{r-1}),$$
which in view of Lemma~\ref{BL} implies
that for $C\sub \Zz^d$ with $|C|<d^{O(1)}$,
$$|\partial C|=\gO(|C|d).$$
Applying this with $C=W\sm G_0$ gives
$|G_0| = (1-O(1/d))g$.  But then
$|\nabla(G_0,A)|\leq \ell|A| =O(|G_0|)$ implies
$$\gd g \ell =|\nabla(G_0,{\cal O}\sm A)|\geq (\ell -O(1))|G_0|
 = \ell (1-O(1/d))g.$$\qed

\subsection{First approximation:  covering the boundary}
\label{CB}

Say a set $C\sub \gG$ {\em separates} $P,Q\sub \gG$ if any path
meeting both $P$ and $Q$ also meets $C$.

In this section we begin the process of approximation by showing
that there is a ``small" collection of subsets of $\gG$, at least
one of which separates $W $ ($=G\cup A$) and $\gG\sm W$ for each
relevant $(G,A)$. We then use these separations to show that there
is a small $\S\sub 2^{\cal E}\times 2^{\cal O}$ such that each of
our $(G,A)$'s is approximated by some $(F,S)\in\S$ in the sense
that \beq{App1.1} S\supseteq A,~F\sub G \enq and \beq{App1.2}
|S\sm A|, |G\sm F|< O(\gd g \sqrt{d\log d}~). \enq This is stated
formally in Lemma~\ref{LApp1} at the end of the section.

Our argument applies to pairs from
$$
\gone :=\{(G,A)\in\g(g,\gd):
\mbox{$(G,A)$ satisfies (\ref{GA0}) and (\ref{G0})}\},
$$
though the main point, Lemma~\ref{ULem}, is valid for all of
$\g(t)$.

In this section (unlike in the next) we make substantial use of
properties particular to $\gG$, specifically the isoperimetric
properties given by Lemma~\ref{BL} and

\beq{quad}
\mbox{$\forall ~w\sim v$ and
$L\sub N(v), ~ |N(w)\cap N(L)|\geq |L|$}
\enq
(which follows from the fact that
for vertices $v\sim w$,
$\gG[(N(v)\cup N(w))\sm \{v,w\}]$ is a matching of
all but one vertex of $N(v)$
and all but one vertex of $N(w)$).

Let
$$G_0' =\{v\in G: d_A(v)\leq \ell/2\} ~~~(\sub G_0),$$
$$B_0' =\{v\in B: d_H(v)\leq \ell/2\} ~~~(\sub B_0),$$
$G_0''=G_0\sm G_0'$ and $B_0''=B_0\sm B_0'$.
Then
\beq{GOBO}
\nabla(G_0'',B_0'')=\0.
\enq
(The more general statement here is:  if $v\in G_0$, $w\in B_0$
and $v\sim w$, then (by (\ref{quad}) with $L=N(v)\cap A$)
$d_G(w )\geq d_A(v) $ ($=\ell-d_B(v)$), implying
$
d_B(v)+d_G(w)\geq \ell.
$)

Notice that (\ref{GOBO}) implies
\beq{sep}
\mbox{$G_0'\cup B_0'$ separates $W$ and $\gG\sm W$}
\enq
(equivalently,
$\nabla(W,\gG\sm W)\sub \nabla(G_0')\cup \nabla(B_0')$).

\begin{lemma}\label{ULem}
In any graph satisfying {\rm (\ref{graph})} and {\rm (\ref{quad})},
for any $(G,A)\in\g(t)$, there exists
$U\sub N(G_0'\cup B_0')$ satisfying

\beq{U1}
N(U)\supseteq \goo\cup\boo
\enq
and
\beq{U2}
|U|<O(t\sqrt{\log \ell/\ell}~).
\enq
\end{lemma}

Before proving this, we observe that it does
accomplish the first goal stated at the
beginning of this section
(existence of a small set of separations).
For $(G,A)$ and $U$ as in Lemma~\ref{ULem},
we have
\beq{U3}
\mbox{$N(U)$ separates $W$ and $\gG\sm W$}
\enq
(by (\ref{sep}) and (\ref{U1})).  So we just need to
limit the number of possibilities for $U$ when $(G,A)\in \gone$.

To do so, notice that \beq{U4} \mbox{$U$ is 6-clustered.} \enq
This follows from Lemma~\ref{Lconn} and (\ref{G0}), once we
observe that $\dist(u,G_0)\leq 2 ~\forall u\in U$ (since $U\sub
N(\goo\cup \boo)$), and that (\ref{U1}) and (\ref{GOBO}) imply
$\dist(v,U)\leq 2~\forall v\in G_0$.

In view of (\ref{U2}) (with $t=\gd g$),
Lemma~\ref{Tree} then gives, for example, a bound
\beq{Ubd}
O(gd^2)(Cd^6)^{O(\gd g \sqrt{\log d/d}~) } =
\exp[O(\gd g d^{-1/2}\log^{3/2}d)]
\enq
on the number of possibilities for $U$.
Here we used Lemma~\ref{deltag} for the equality in
(\ref{Ubd}).
The initial $O(gd^2)$ corresponds to a choice of $x_0$ in
Lemma~\ref{Tree}:
in view of (\ref{GA0}),
there must be some $j\in\indices$ and $k\leq g/(2d)$ for which
$y_0:=v_0+(2k-1)e_j\in G_0$; there are at most $g$ possibilities
for this $y_0$, so at most $O(gd^2)$ possibilities
for a vertex $x_0$ with $d(x_0,y_0)\leq 2$;
and by (\ref{U1}) and (\ref{GOBO})
$U$ must contain such an $x_0$.

\mn
{\em Proof of Lemma~\ref{ULem}.}

By ``duality" (see Section~\ref{Pre})
it's enough to show the existence of $S\sub N(\goo)$
with
\beq{S1}
N(S)\supseteq \goo
\enq
and
\beq{S2}
|S|< O(t\sqrt{\log \ell/\ell}~).
\enq

Define
$Q=\{v\in G_0:d_A(v)\leq \sqrt{\ell\log \ell} \}$,
$K=G_0\sm Q$,
and $P=N(Q)\cap A$.
By (\ref{quad}),
\beq{degP}
d_{G_0}(v)\geq \ell-\sqrt{\ell\log \ell} ~~\forall v\in P.
\enq
Let $P'=\{v\in P:d_K(v)\geq \ell/2\}$, $P''=P\sm P'$,
$Q'=Q\cap N(P')$, $Q''=Q\sm Q'$ and
$R=\{v\in B_0\cap N(G_0'):d_{G_0}(v)> \sqrt{\ell\log \ell} \}$.

Now $P''$ is a cover of $Q''$ of size
$O(t\sqrt{\log \ell/\ell}~)$, the size bound following from
$|Q|\leq t\ell/(\ell- \sqrt{\ell\log \ell}~) =O(t)$
(using (\ref{nabla})),
$d_{P''}(v)\leq d_A(v)\leq \sqrt{\ell\log \ell} ~~\forall v\in Q$,
and $d_Q(v)>\ell/2-\sqrt{\ell\log \ell}  ~~\forall v\in P''$
(using (\ref{degP}) and the definition of $P''$).

On the other hand, we can cover $\goo \sm Q''$ by a similarly
small subset of $R$, as follows.
From (\ref{quad}) we have $N(K)
\cap N(G_0')
\cap B_0\sub R$.
This gives
$d_R(v)>\ell/2$ for $v\in \goo\sm Q$, while for $v\in Q'$,
$$
d_R(v)\geq |N(v)\cap N(K)| -|N(v)\cap A|\geq \ell/2
-\sqrt{\ell\log \ell}
$$
(the second inequality following from (\ref{quad}) and the
definitions of $Q'$ and $Q$).
So, noting that $|R|<t\sqrt{\ell/\log \ell}$
(again using (\ref{nabla})),
Lemma~\ref{Lcor} says that we can cover $G_0'\sm Q''$ by
some $T\sub R$ of size at most
$|R|(1+\log \ell) /(\ell/2  -\sqrt{\ell\log \ell}~)
<O(t\sqrt{\log \ell/\ell}~)$.
(And note $P\sub N(G_0')$ since $Q\sub G_0'$, and
$R \sub N(G_0')$ by definition, so $S:=P''\cup T\sub N(G_0')$.)
\qed

We now return to $\gG$.
Given $U$ as above, let us temporarily set $L=N(U)$.
Then
$|L| =O(\gd g \sqrt{d\log d}~)$.

Say a component $C$ of $\gG-L$ is {\em large} if $|C|>d$ and {\em
small} otherwise. Lemma~\ref{BL} implies

$$|\nabla(C,L)| =|\nabla(C)| \geq |\partial C|=  \gO(|C|d)$$
for small $C$ (actually also for considerably larger $C$), and

$$|\nabla(C,L)|= \gO(d^2)$$
for large $C$.
But $|\nabla(L)|\leq 2d|L|= O(\gd g d^{3/2}\sqrt{\log d}~ )$,
so
\beq{largenumber}
\mbox{the number of large components is
$O(\gd g d^{-1/2}\sqrt{\log d}~  )$,}
\enq
and the number of vertices in small components is
$O(\gd g \sqrt{d\log d}~  )$.

It follows that if $(G,A)$ is any pair satisfying (\ref{GA3})
for which $L$ separates
$W  $ and $\gG\sm W$, then we
satisfy (\ref{App1.1})
and (\ref{App1.2}) with

\beq{FS1} \mbox{$F=P\cap {\cal E}~~$ and $~~  S= (P\cup Q\cup
L)\cap {\cal O}$,} \enq where $P$ is the union of those large
components of $\gG-L$ that meet (equivalently, are contained in)
$W$, and $Q$ is the union of (all) the small components. In
particular this is true if $(G,A)$ is any pair from $\gone$ for
which Lemma~\ref{ULem} applied to $(G,A)$ produces $U$.

By (\ref{largenumber}) the number of possibilities (given $L$)
for $(F,S)$ as in (\ref{FS1}) is at most
$\exp[O(\gd g d^{-1/2}\sqrt{\log d~} )]$, and combining this
with the bound (\ref{Ubd}) on the number of $U$'s we have

\begin{lemma}\label{LApp1}
There exist $\S\sub 2^{\cal E}\times 2^{\cal O}$ with \beq{Sbd}
|\S|<\exp[O(\gd g d^{-1/2}\log^{3/2}d)] \enq and a map
$\pi_1:\gone\ra\S$ such that {\rm (\ref{App1.1})} and {\rm
(\ref{App1.2})} hold for each $(G,A)\in \gone$ and
$(F,S)=\pi_1(G,A)$.
\end{lemma}

\subsection{Second approximation}
\label{SA}

The discussion in this section
is valid for any graph $\gS$ satisfying (\ref{graph}).
It may be worth reiterating that we follow the conventions
given at the end of Section~\ref{Pre}.

Given $(F^*,S^*)\in 2^{\cal E}\times 2^{\cal O}$ and a positive
$x$, write $\gtwo=\gtwo(F^*,S^*,x)$ for the set of $(G,A)$'s in $
\g(t)$ satisfying (\ref{App1.1}) (with $(F^*,S^*)$ in place of
$(F,S)$) and \beq{SF2} |S^*\sm A|, |G\sm F^*|<x. \enq

\begin{lemma}\label{App2}
With notation as above, for any $0<\psi<\ell$, there exist $\T\sub
2^{\cal E}\times 2^{\cal O}$, \beq{T} |\T| <
\exp[O((x/\ell)+(t/\psi)) \log \ell], \enq and a map
$\pi_2:\gtwo\ra\T$ such that for each $(G,A)\in\gtwo$ and
$(F,S)=\pi_2(G,A)$ we have {\rm (\ref{App1.1})} and \beq{App2.2}
v\in S\Rightarrow d_{F}(v)>\ell-\psi, ~~~v\in E\Rightarrow
d_{T}(v)>\ell-\psi \enq (where as usual $E={\cal E}\sm F$ and
$T={\cal O}\sm S$).
\end{lemma}

\mn {\em Remarks.} We only need Lemma~\ref{App2} when
$(F^*,S^*)\in \S$ (with $\S$ as in Lemma~\ref{LApp1}), in which
case we take $t=\gd g$ and $x=O(\gd g \sqrt{d\log d}~)$ (with an
appropriate constant), so that $\gtwo\supseteq
\pi_1^{-1}(F^*,S^*)$; but the extra generality costs us nothing.
The pairs we produce will satisfy $S\sub S^*$ and $F\supseteq
F^*$, but we don't need this in what follows.


\mn
{\em Proof of Lemma~\ref{App2}}

We would like to exhibit a procedure which, for a
given $(G,A)\in \gtwo$, outputs a pair $(F,S)$ satisfying
(\ref{App1.1}) and (\ref{App2.2}),
and show that the set $\T$ of pairs produced in this way
is small.

We produce $(F,S)$ via a sequence of modifications,
initializing at $(F,S)=(F^*,S^*)$.
Note that whenever we update $(F,S)$,
we also automatically update $E,T$, etc.

One preliminary observation: \beq{preob} |S_0^*|, |E_0^*|<  x
+\ell x \enq (since $S_0^* \sub (S^*\sm A)\cup N(G\sm F^*)$, and
similarly for $E_0^*$; recall $S_0^*=S^* \cap N(E^*)$ and
$E_0^*=E^* \cap N(S^*)$, where $E^*={\cal E}\sm F^*$).

\mn
{\bf Stage 1A}
Set $\xi = \ell/2$.

\mn
{\bf (A.1)}
Repeat for as long as possible:
choose $w\in H$ with $d_S(w)\geq\xi$ and
do $S\leftarrow S\sm N(w)$.

\mn {\bf (A.2)} When no longer possible, do $F\la F\cup\{w\in
{\cal E}: d_S(w) \geq\xi\}$.

\mn
{\bf Stage 1B}
Do the same thing in the dual; that is,

\mn
{\bf (B.1)}
for as long as possible, choose
$w\in A$ with $d_E(w)\geq\xi$
and do $F\la F\cup N(w)$, and

\mn {\bf (B.2)} when no longer possible, do $S\la S\sm \{w\in
{\cal O}:d_E(w)\geq\xi\}$.

\mn
Notice---a crucial idea---that $(F,S)$ produced by
Stage 1 does satisfy (\ref{App1.1}).


\mn
{\bf Analysis:}

The output $(F,S)$ of Stage 1 is determined by
the sets of $w$'s used in (A.1) and (B.1).

Since each iteration in (A.1)
shrinks $|S|$ by at least  $\xi $ while
maintaining $A\sub S$, the number of iterations
is less than $x/\xi =2x/\ell$.
Moreover, each $w$ used in (A.1) lies in
$N(S_0^*)$.  So the number of possibilities for the
set of $w$'s used in (A.1) is less than
$\sum_{i\leq x/\xi}\C{\ell|S_0^*|}{i} <
\exp[O((x/\ell)\log \ell)]$
(using (\ref{preob})).

At the end of (A.2) we have
$w\in G\sm F \Ra d_T(w)>\ell-\xi =\ell/2$,
which, since $| \nabla(G,T)| \leq t\ell$ (see (\ref{nabla})), gives
$|G\sm F| < 2t$.

Similarly, the number of choices for the set of
$w$'s used in Stage 1B is at most $\exp[O((x/\ell)\log \ell)]$
(note Stage 1A does not increase $E_0^*$),
and at the end of this stage we have
$|S\sm A|  <2t$.

\bigskip
Stage 2 now repeats Stage 1, starting with
the revised $(F,S)$, using $\psi $
in place of $\xi$, and replacing (\ref{SF2})
and (\ref{preob}) by
$$|S\sm A|, |G\sm F|<2t$$
and
$$|S_0|, |E_0|<  2t(1+\ell).$$

This clearly produces an $(F,S)$ satisfying (\ref{App1.1}) and
(\ref{App2.2}). Moreover, repeating the analysis above, we find
that the number of possible outputs of Stage 2, for a given output
of Stage 1, is at most $\exp[O((t/\psi)\log \ell)]$. So the number
of possible outputs of the entire procedure is no more than $
\exp[O((x/\ell)+(t/\psi)) \log \ell] $.\qed

\subsection{Status}
\label{Status}

We now specify $t=\gd g$ and $x=O(\gd g \sqrt{d\log d}~)$
(the bound in (\ref{App1.2})), and
$\psi =\sqrt{d}$ (any
$\psi \in (\gO(\sqrt{d/\log d}), O(\sqrt{d\log d}))$
would do; see the remark following (\ref{nuJlast}).)
Specializing to these values and combining
Lemmas~\ref{LApp1} and \ref{App2}, we have

\begin{lemma}\label{LApp3}
There exist $\U\sub 2^{\cal E}\times 2^{\cal O}$, \beq{Ubd2}
|\U|<\exp[O(\gd g d^{-1/2}\log^{3/2}d)], \enq and $\pi:\gone\ra\U$
such that {\rm (\ref{App1.1})} and {\rm (\ref{App2.2})} hold for
each $(G,A)$ and $(F,S)=\pi(G,A)$.
\end{lemma}
(The expression in the exponent in (\ref{Ubd2}) is
the maximum of
the corresponding expressions from (\ref{Sbd}) and (\ref{T}).)

Now consider some $(F,S)\in\U$.
Notice that, for any $(G,A)\in\pi^{-1}(F,S)$,
$Q:=S_0\cup E_0$ contains all
vertices whose locations in the partition
$\gG=G\cup H \cup A\cup B$
are as yet unknown;
namely, we have
$$
F\sub G, ~~ T\sub B, ~~ S\sm S_0\sub A,~~ E\sm E_0\sub H
$$
(the first two containments
are just (\ref{App1.1});
$S\sm S_0\sub A$ follows from $F\sub G$,
(\ref{GA3}) and the definition
of $S_0$, and $ E\sm E_0\sub H$ is similar).

{\em By convention, whenever we are given an
$(F,S)$, we take $Q$ to
be as defined in the preceding paragraph, and
write
$\gG_Q$ for the subgraph induced by $Q$.}

\subsection{Flow}
\label{Flow}

Here, finally, we define $\nu$
(for large $I$; for small $I$,
see Section \ref{Pre}).

Throughout the section we fix $(F,S)\in\U$.
It is now convenient to write $G\sim (F,S)$ if $\pi(G,A)=(F,S)$
and $I\sim (F,S)$ if $G(I)\sim (F,S)$.

To define $\nu (I,\cdot)$ for
$I\sim (F,S)$,
we first need to choose a direction $j=j(I)$.
Fix such an $I$ and
let $G=G(I)$, $A=A(I)$, etc.  The choice of $j$ will
depend only on $(G,A)$.
Observe that (using (\ref{App2.2}))
$$
\sum_j|\gs_j(S_0\cap A)\cap E_0|=|\nabla(S_0\cap A, G\cap E_0)|
<|G\cap E_0|\psi,
$$
$$
\sum_j |\gs_j^{-1}(E_0)\cap (S_0\sm A)| = |\nabla(E_0,S_0\sm A)|
< |S_0\sm A|\psi.
$$
But (\ref{App2.2}) and (\ref{nabla}) imply
$ |G\cap E_0|+|S_0\sm A|< \gd g \ell/(\ell-\psi) $,
so that
\begin{eqnarray}\label{junk}
\sum_j |\gs_j(S_0)\cap E_0| &=&
\sum_j
(|\gs_j(S_0\cap A)\cap E_0|+|\gs_j^{-1}(E_0)\cap (S_0\sm A)|)
\nonumber\\
&<& \gd g \ell\psi/(\ell-\psi).
\end{eqnarray}

We assert that we can choose $j$ so that
\beq{j1}
|G_0^j|>.8\gd g
\enq
and
\beq{j2}
|\gs_j(S_0)\cap E_0|
<10 |G_0^j|\psi/\ell.
\enq
To see this, let
$$
P=\{j\in [-d,d]\sm\{0\}:
|\gs_j(S_0)\cap E_0|
\geq 10 |G_0^j|\psi/\ell\}.
$$
Then (\ref{junk}) gives
$$
\sum_{j\in P} |\goj|\leq \frac{\ell}{10\psi}\sum
|\gs_j(S_0)\cap E_0|< \gd g \frac{\ell^2}{10(\ell-\psi)},
$$
so (using (\ref{Gsum}))
$$\sum_{j\not\in P}|\goj|>(1-\ell/(10(\ell-\psi)))\gd g \ell.$$
So there exists $j\not\in P$ with (say) $|\goj|>.8\gd g$,
which is what we want.\qed

Having chosen $j$ satisfying (\ref{j1}) and (\ref{j2}),
we turn to defining $\nu(I,\cdot)$.  Let
$$
C = C^j(I) =\goj\cap F\cap \sj(S_0) ~~(=\sj(S_0\sm A)\cap F),
$$
\begin{eqnarray*}
D=D^j(I) &=& \goj\cap (\sj(T)\cup (\sj(S_0)\cap E_0)).
\end{eqnarray*}
Then
\beq{CD}
\mbox{$C\cup D$ is a partition of $\goj$.}
\enq

Setting
$\ga=\ga(\gl)=\gl /(1+\gl)^2$ and
$\gb=\gb(\gl) = 1-\ga\gl =(1+2\gl)/(1+\gl)^2$,
define
$$
\nu(I,J)=\left\{\begin{array}{ll}
(\ga\gl)^{|C\cap J|}
\gb^{|C\sm J|}(\gl/(1+\gl))^{|D\cap J|}(1+\gl)^{-|D\sm J|}&\\
~~~~~~~~~~~=\frac{w(J)}{w(I)}
\ga^{|C\cap J|}
\gb^{|C\sm J|}
(1+\gl)^{-|D|}&\mbox{if $j\in\varphi_j(I)$}\\
0&\mbox{otherwise.}
\end{array}\right.
$$
Then
\beq{Isum}
\sum_J\nu(I,J) =1 ~~~~~\forall I
\enq
(because of (\ref{CD})).
On the other hand we will show, for any $J$,
\beq{MAIN}
\sum_{I\sim (F,S)}\frac{w(I)}{w(J)}\nu(I,J) \leq \ell\gb^{ \gd g/2}.
\enq

\subsection{Proof of (\ref{MAIN})}
\label{POM}

We need one easy lemma.  Given a bigraph $\Sigma$ on $P\cup R$ and
$U\sub R$, say that a (vertex) cover $K\cup L\cup M$ of $\Sigma$
with $K\sub P$, $L\sub U$ and $M\sub R\sm U$ is {\em legal}
(with respect to $U$) if it is a minimal cover and
$$K=N(U\sm L).$$
(Note minimality implies $K=N(R\sm (L\cup M))$.)

\begin{lemma}\label{Llegal}
With notation as above, let $K\cup L\cup M$ be a legal cover with $|K\cup L|$
as small as possible.  Then

\mn
{\rm (a)}  $\forall K'\sub K  ~~~~|N(K')\cap (U\sm L)|\geq |K'|$,

\mn
{\rm (b)}  $\forall L'\sub L  ~~~~|N(L')\sm K|\geq |L'|$.
\end{lemma}
{\em Proof.}
(a)  Given $K'\sub K$,
let $S=N(K')\cap (U\sm L)$,
$$K'' =\{v\in K: N(v)\cap U\sub S\cup L\} ~~(\supseteq K'),$$
and $T=N(K'')\cap (R\sm U)$.
Then

\mn
(i)  $(K\sm K'')\cup (L\cup S)\cup(M\cup T)$ is a minimal cover

\mn
(a straightforward verification using the fact that each vertex of
$K\sm K''$ has a neighbor in $U\sm (L\cup S)$), and

\mn
(ii) $K\sm K''=N(U\sm (L\cup S))$.

\mn
Minimality of $|K\cup L|$ thus implies
$|K\sm K''|+|L\cup S|\geq |K|+|L|$,
so $|S|\geq |K''|\geq |K'|$.

\mn
(b)  This is similar.
Given $L'\sub L$, let
$W=N(L')\sm K $ and
$$L'' =\{u\in L\cup M: N(u)\sub K\cup W\} ~~(\supseteq L').$$
Then

\mn
(i)  $K\cup W\cup ((L\cup M)\sm L'')$
is a minimal cover, and

\mn
(ii) $K\cup W =N(U\sm (L\sm L''))$.

\mn
Minimality of $|K\cup L|$ thus implies
$|K\cup W|+|L\sm L''|\geq |K|+|L|$,
and $|W|\geq |L''|\geq |L'|$.\qed

\bn
{\em Proof of} (\ref{MAIN}).

Given $(F,S)$, $J$ and $j$, set
$$\ist=\ist(F,S,J,j) =\{I\sim (F,S): j(I)=j, J\in\varphi_j(I)\}.$$
We will show
$$
\sum_{I\in\ist}\frac{w(I)}{w(J)}\nu(I,J)<\gb^{ \gd g/2},
$$
which of course gives (\ref{MAIN}).

Set $U=\sji(J)\cap S_0$.
Suppose $I\in \ist$, and set $G=G(I)$, $A=A(I)$, and
$$
K=K(I)=G\cap E_0,
~~~L=L(I) =U\sm A, ~~
~ M =M(I) =(S_0\sm U)\sm A.
$$
Then $K\cup L\cup M$ ($ = (G\cup B)\cap Q$)
is a minimal cover of $\gcq$.
(That it is a cover follows from (\ref{GA3}); for minimality,
notice (e.g.) that each $v\in G\cap E_0$ has a neighbor in $A$,
which must be in $S_0$ (using $A\sub S$ and the definition of
$S_0$).)
Moreover, we assert,
\beq{KNQ}
K=N_{\Gamma_Q}(U\sm L).
\enq

\mn
{\em Proof.}
We show that each side of (\ref{KNQ}) contains the other.
The obvious direction is
$$
N_{\Gamma_Q}(U\sm L)=N_{\Gamma_Q}(U\cap A)\sub N(A)\cap E_0 =
G\cap E_0 =K.
$$
For the reverse containment, suppose $v\in K$.
Since $K\sub G_0$,
(\ref{GA7}) says that
$v$ has a neighbor
$u\in A\cap I$.
Then $u\in S_0$ (because $v\in E_0\not\sim S\sm S_0$),
implying $u\in U$ (since $u\in A\cap I\Rightarrow \sj(u)\in J$).
And of course $u\not\in L$ (since $u\in A$).\qed

Thus $K\cup L\cup M$ is a legal cover of $\gcq$ with respect to $U$
in the sense of Lemma~\ref{Llegal}.

Now fix $K_0\cup L_0\cup M_0$, a legal cover of $\gcq$ with respect to $U$
with $|K_0\cup L_0|$ as small as possible.

Given $I\in\ist$, let $K=K(I)$ etc. be as above and set
$K'=K_0\sm K$,
$L'=L_0\sm L$.
Then by Lemma~\ref{Llegal},
\beq{LK}
|L|\geq |K'|+|L_0\sm L'|,  ~~~
|K|\geq |L'|+|K_0\sm K'|.
\enq

Furthermore, we assert,
\beq{knowK}
K=(K_0\sm K')\cup N_{\gcq}(L').
\enq
The point of this is that it says that $(K',L')$ determines
$G$ (so also $A$), and therefore $I\in \ist$
(because of (\ref{recover1})).

To see (\ref{knowK}), just observe that the only point
requiring proof is $K\sm K_0\sub N_{\gcq}(L_0\sm L)$,
and that this follows from (\ref{KNQ}) once we notice that
$\nabla(K\sm K_0,U\sm (L_0\cup L))=\0$
(since $K_0\cup L_0$ covers $\nabla(E_0, U)$).

Now with $C=C^j(I)$, $D=D^j(I)$ as in the discussion
preceding (\ref{CD}), observe that
$$\mbox{$C\cap J=\sj(L\sm \sji(E_0))~~~$ and
$~~~C\sm J=\sj(M\sm \sji(E_0))$,}$$
and that we may partition $D$ as
$$
D=(\sj(T)\cap F)\cup (K\sm \sj(S_0\sm (L\cup M))).
$$
Thus, with inequalities justified below,
\begin{eqnarray}
\frac{w(I)}{w(J)}\nu(I,J) &=&
\ga^{|\sj(L\sm \sji(E_0))|}\gb^{|\sj(M\sm \sji(E_0))|}\nonumber\\
&& ~~~~~~\cdot
(1+\gl)^{-(|\sj(T)\cap F|+ |K\sm \sj(S_0\sm (L\cup M))|)}
\nonumber\\
&\leq& \ga^{|L|}\gb^{|M|}(1+\gl)^{-(|K|+|\sj(T)\cap F|)}\nonumber\\
&& ~~~~~~\cdot
\ga^{-(|\sj(S_0\cap A)\cap K|+|\sji(E_0)\cap (S_0\sm A)|)}
\label{Jsum1}\\
&\leq &
\ga^{|L|}(1+\gl)^{-|K|}\gb^{|\goj|-(|K|+|L|)}
\ga^{-O(|\goj|\psi/\ell)}\label{Jsum2}\\
&\leq &
\gb^{\gd g/2}\ga^{|L|}(1+\gl)^{-|K|}\gb^{-(|K|+|L|)}
\label{Jsum3}\\
&=&\gb^{\gd g/2} \left(\frac{1+\gl}{1+2\gl}\right)^{|K|}
\left(\frac{\gl}{1+2\gl}\right)^{|L|}
\nonumber\\
&\leq &
\gb^{\gd g/2} \left(\frac{1+\gl}{1+2\gl}\right)^{|L'|+|K_0\sm K'|}
\left(\frac{\gl}{1+2\gl}\right)^{|K'|+|L_0\sm L'|}
\label{Jsum4}\\
&=&
\gb^{\gd g/2} \left(\frac{1+\gl}{1+2\gl}\right)^{|K_0|}
\left(\frac{\gl}{1+2\gl}\right)^{|L_0|}
\left(\frac{\gl}{1+\gl}\right)^{|K'|-| L'|}.\nonumber
\end{eqnarray}
(In
(\ref{Jsum1}) we used
$\ga^{-1}=\max\{\ga^{-1},\gb^{-1},1+\gl\}$;
in
(\ref{Jsum2}) we used
$\goj\sub\sj(L\cup M)\cup K\cup (\sj(T)\cap F)$,
$(1+\gl)^{-1}<\gb$ and (\ref{j2});
(\ref{Jsum3}) is from (\ref{j1}), using
$(\psi/\ell)\log(1/\ga)=o(\log(1/\gb))$, which is a consequence of

\beq{lampsi}
\gl^2 =\omega((\psi/\ell)\log(1/\gl))
\enq
for small $\gl$, and easily verified when $\gl$ is larger;
and
(\ref{Jsum4}) comes from (\ref{LK}).)

Thus, recalling---see the remark
following (\ref{knowK})---that
each $(K',L')$ corresponds to at most one $I\in \ist$,

\begin{eqnarray*}
\sum_{I\in\ist}\frac{w(I)}{w(J)}\nu(I,J) &\leq &
\gb^{\gd g/2} \left(\frac{1+\gl}{1+2\gl}\right)^{|K_0|}
\left(\frac{\gl}{1+2\gl}\right)^{|L_0|}
\sum_{K'\sub K_0}
\sum_{L'\sub L_0}\left(\frac{\gl}{1+\gl}\right)^{|K'|-| L'|}\\
&=&\gb^{\gd g/2}.
\end{eqnarray*}

As noted earlier this gives (\ref{MAIN}).

\subsection{Finally}
\label{Fin}

Now fixing $J\in\J$, we are ready to verify
(\ref{nuin}) (thus completing the proofs
of Theorems~\ref{Thm} and \ref{Thm0}).

Note first of all (referring to (\ref{Ubd2})) that for $\gl \leq
2$ (say) (\ref{MAIN}) implies

\begin{eqnarray}
\sum_{I\in\J(g,\gd)}\frac{w(I)}{w(J)}\nu(I,J) &=& \sum_{(F,S)\in
\U}\sum_{I\sim(F,S)} \frac{w(I)}{w(J)}\nu(I,J)
\nonumber\\
&\leq& |\U|\ell\gb ^{\gd g/2}
\nonumber\\
&<&\ell\exp[\{O(d^{-1/2}\log^{3/2}d) -\Omega(\gl^2)\}\gd g]
\nonumber\\
&<&\exp[-\Omega(\gl^2\gd g)],\label{nuJlast}
\end{eqnarray}
while for larger $\gl$,
\begin{equation} \label{largelambda}
\sum_{I\in\J(g,\gd)}\frac{w(I)}{w(J)}\nu(I,J) <\gl^{-\gO(\gd g)}.
\end{equation}

\mn {\em Remark.} Our choice of $\psi$ was constrained by the
demands of (\ref{lampsi}) and (\ref{nuJlast}) (the latter since
$\psi = o(\sqrt{d/\log d})$ would give---via (\ref{T})---a larger
bound in (\ref{Ubd2})).

We first deal with large $I$'s (recall $I$ is large if
$|G(I)|>d^3$). Here we have already done the work: Assuming first
that $\gl\leq 2$, and with justifications to follow, we have

\begin{eqnarray}
\sum_{\mbox{{\small $I$ large}}}
\frac{w(I)}{w(J)}\nu(I,J)
& = &
    \sum_{g> d^3}\sum_{\gd}\sum_{I\in I(g,\delta)}
\frac{w(I)}{w(J)}\nu(I,J) \nonumber \\
 & = &     \sum_{g> d^3}\sum_{\gd}
\exp[-\gO(\gl^2 \delta g)] \label{bylast}\\
 & \leq & \sum_{g> d^3}\sum\{\exp[-\gO(\gl^2 i)]:
i\geq \gO(d^{-1}g^{1-1/d})\} \label{delbig} \\
 & \leq & \sum_{g> d^3}
\exp[-\gO(\gl^2(d^{-1}g^{1-1/d}))] \label{xl} \\
 &<&
\exp[-\gO(\gl^2d^{3(1-1/d)-1})] \label{sums} \\
&<&
\exp[-\go(\gl d)].\label{note3}
\end{eqnarray}
Of course sums involving $\gd$, are
restricted to $\gd$ for which $\gd g$ is
an integer.
The main inequality
(\ref{bylast})
is just (\ref{nuJlast}), and
(\ref{delbig}) comes from
Lemma~\ref{deltag}.
In (\ref{xl}) we have absorbed a
factor $\gl^{-2}$ in the exponent.
One way (probably not the most natural) to see the
inequality in (\ref{sums}) is to use
$$
\mbox{$(1-\eps)^{g^{1-\gd}} < (1-\eps)^{iK^{1-\gd}} ~~$
for $~ i^{1/(1-\gd)}K < g\leq (i+1)^{1/(1-\gd)}K$}
$$
with $K=d^3$, $\gd = 1/d$ and $1-\eps = \exp[-\gO(\gl^2d^{-1})]$.

For $\gl >2$ a similar analysis (using (\ref{largelambda})) gives
\beq{68'} \sum_{\mbox{{\small $I$ large}}}
\frac{w(I)}{w(J)}\nu(I,J) \leq \gl^{-\gO(d^2)}. \enq

\medskip
Finally we turn to the easy case of small $I$.
Here we abuse our notation slightly and set
$$\J(g,a) =\{I\in\J_0:|G(I)|=g,|A(I)|=a\}.$$
For a (nonempty) $\J(g,a)$ with $g<d^3$, Lemma~\ref{deltag} gives
$a=O(g/d)$, so that, since each $A(I)$ is 2-clustered and contains
$v_0$, Lemma~\ref{Tree} bounds the number of possibilities for
$A(I)$ with $I\in\J(g,a)$ by $\exp[O((g/d)\log d)]$.

But we also know (see (\ref{recover1})) that, given $J$ and $j$,
$I\in\varphi_j^{-1}(J)$ is determined by $G(I)$ (or $A(I)$),
and that (by (\ref{easynu}), (\ref{gojbig}), and again
Lemma~\ref{deltag})
\begin{eqnarray*}
\frac{w(I)}{w(J)}\nu(I,J)&=& (1+\gl)^{-|\goj(I)|} \\
&\leq &(1+\gl)^{-\gd g}\\
&=&(1+\gl)^{-(1-O(1/d))g}.
\end{eqnarray*}

So finally, noting that $A(I)\neq\0$ implies $|G(I)|\geq \ell$,
we have
\begin{eqnarray*}
\sum_{I\in\J(g,a)}
\frac{w(I)}{w(J)}\nu(I,J) &<&
 \ell\exp[O((g/d)\log d)] (1+\gl)^{-(1-O(1/d))g} \\
  &<& (1+\gl)^{-(1-o(1))g}
\end{eqnarray*}

and

\begin{eqnarray*}
\sum_{\mbox{{\small $I$ small}}}
\frac{w(I)}{w(J)}\nu(I,J) & = &
     \sum_{\ell\leq g\leq d^3}\sum_{a\leq g}\sum_{I\in\J(g,a)} \frac{w(I)}{w(J)}\nu(I,J)
\nonumber \\
 & < & \sum_{\ell\leq g<d^3}g (1+\gl)^{-(1-o(1))g}
\nonumber \\
 & \leq & (1+\gl)^{-(1-o(1))\ell};
\end{eqnarray*}
and combining this with (\ref{note3}) or (\ref{68'}) gives (\ref{nuin}).

\bn
{\bf Acknowledgments}
We are very grateful to Vladimir Gurvich for translating
parts of \cite{Sapozhenko}.  Thanks also
to Chuck Weibel for help
with the proof of Lemma~\ref{connbound}.

\end{document}